\documentclass[titlepage,12pt]{article}
\usepackage{amssymb}
\usepackage{amsfonts}
\begin{document}


{\Large \bf Rigid and Non-Rigid Mathematical \\ \\
            Theories : the Ring $\mathbb{Z}$ Is Nearly Rigid} \\

{\bf Elem\'{e}r E Rosinger} \\
Department of Mathematics \\
and Applied Mathematics \\
University of Pretoria \\
Pretoria \\
0002 South Africa \\
eerosinger@hotmail.com \\

\hfill {\it Dedicated to Marie-Louise Nykamp} \\ \\

{\bf Abstract} \\

Mathematical theories are classified in two distinct classes : {\it rigid}, and on the other hand, {\it non-rigid} ones. Rigid theories, like group theory, topology, category theory, etc., have a basic concept - given for instance by a set of axioms - from which all the other concepts are defined in a unique way. Non-rigid theories, like ring theory, certain general enough pseudo-topologies, etc., have a number of their concepts defined in a more free or relatively independent manner of one another, namely, with {\it compatibility} conditions between them only. As an example, it is shown that the usual ring structure on the integers $\mathbb{Z}$ is not rigid, however, it is nearly rigid. \\ \\

{\large \bf 0. Introduction} \\

Rigid theories, like group theory, topology, category theory, etc., have a basic concept - given for instance by a set of axioms - from which all the other concepts are defined in a unique way. Non-rigid theories, like ring theory, certain general enough pseudo-topologies, etc., have a number of their concepts defined in a more free or relatively independent manner of one another, namely, with {\it compatibility} conditions between them only. \\

One can note that even in Algebra there are nonrigid mathematical structures. For instance, let $( R, +, . )$ be a ring. Then in principle, neither the addition "+" determines the multiplication ".", nor multiplication determines addition. Instead, they are relatively independent of one another, and only satisfy the usual compatibility
conditions, namely, the distributivity of multiplication with respect to addition. \\

On the contrary, in groups $( G, \diamond )$, all concepts are defined uniquely based eventually on the underlying set $G$ and the binary operation $\diamond$. \\

Non-rigid mathematical structures need not always form usual Eilenberg - Mac Lane categories, [8,12], but more general ones, as illustrated by the case of certain general enough concepts of pseudo-topology. \\

Such a rather general concept of pseudo-topology was used in constructing differential algebras of generalized functions containing the Schwartz distributions, [1-7,10]. These algebras proved to be convenient in solving large classes of nonlinear partial differential equations, see [10] and the literature cited there, as well as section 46F30 in the Subject Classification 2009 of the American Mathematical Society, at
www.ams.org/msc/46Fxx.html \\

And it is precisely because of that non-rigid character that the totality of such pseudo-topologies does {\it no longer} constitute a usual Eilenberg - Mac Lane category, but one which is more general, [8,12]. \\

As it happens, the rigid structure of the usual Hausdorff-Kuratowski-Bourbaki, or in short, HKB concept of topology is also one of the reasons for a number of its important deficiencies, such as for instance, that the category of such topological spaces is not Cartesian closed. \\

Spaces $( \Omega, {\cal M}, \mu )$ with measure, where $\Omega$ is the underlying set, ${\cal M}$ is a $\sigma$-algebra on it, and $\mu : {\cal M} \longrightarrow \mathbb{R}$ is a $\sigma$-additive measure, are further examples of non-rigid structures, since for a given $( \Omega, {\cal M} )$, there can in general be infinitely many associated $\mu$. \\

Topological groups, or even topological vector spaces, are typically non-rigid structures. Indeed, on an arbitrary group, or even vector space, there may in general be many compatible topologies, and even Hausdorff topologies. \\

Obviously, an important advantage of a rigid mathematical structure, and in particular, of the usual HKB concept of topology, is a simplicity of the respective theoretical development. Such simplicity comes from the fact that one can start with only one single concept, like for instance the open sets in the case of HKB topologies, and then based on that concept, all the other concepts can be defined in a unique manner. \\
Consequently, the impression may be created that one has managed to develop a universal theory in the respective discipline, universal in the sense that there may not be any need for alternative theories in that discipline, as for instance is often the perception about the HKB topology. \\

The disadvantage of a rigid mathematical structure is in a consequent built in lack of flexibility regarding the interdependence of the various concepts involved, since each of them, except for a single starting concept, are determined uniquely in terms of that latter one. And in the case of the HKB topologies this is manifested, among others, in the difficulties related to dealing with suitable topologies on spaces of continuous functions, that is, in the failure of the category of such topological spaces to be Cartesian closed. \\

Non-rigid mathematical structures, and in particular, certain general enough pseudo-topologies, can manifest fewer difficulties coming from a lack of flexibility. \\

A disadvantage of such non-rigid mathematical structures - as for instance with various approaches to pseudo-topologies - is in the large variety of ways the respective theories can be set up. Also, their respective theoretical development may turn out to be more complex than is the case with rigid mathematical structures. \\
Such facts can lead to the impression that one could not expect to find a universal enough non-rigid mathematical structure in some given discipline, and for instance, certainly not in the realms of pseudo-topologies. \\

As it happens so far in the literature on pseudo-topologies, there seems not to be a wider and explicit enough awareness about the following two facts

\begin{itemize}

\item one should rather use non-rigid structures in order to avoid the difficulties coming from the lack of flexibility of the rigid concept of usual HKB topology,

\item the likely consequence of using non-rigid structures is the lack of a sufficiently universal concept of pseudo-topology.

\end{itemize}

As it happens, such a lack of awareness leads to a tendency to develop more and more general concepts of pseudo-topology, hoping to reach a sufficiently universal one, thus being able to replace once and for all the usual HKB topology with "THE" one and only "winning" concept of pseudo-topology. \\
Such an unchecked search for increased generality, however, may easily lead to rather meagre theories. \\

It also happens in the literature that, even if mainly intuitively, when setting up various concepts of pseudo-topology a certain restraint is manifested when going away from a rigid theory towards some non-rigid ones. And certainly, the reason for such a restraint is that one would like to hold to the advantage of rigid theories which are more simple to develop than the non-rigid ones. \\

Amusingly, the precedent in Geometry, happened two centuries earlier, is missed from the view both by those who hold to the usual HKB concept of topology, as well as by those trying to set up a general enough pseudo-topology which hopefully may be universal. After all, having by now gotten accustomed that, in fact so fortunately, Geometry can mean many things in different situations, it may perhaps be appropriate to accept a similar view regarding Topology ... \\

Lastly, let us note that in modern Mathematics it is "axiomatic" that theories are built as {\it axiomatic
systems}. \\
As it happens, however, ever since the early 1930s and G\"{o}del's Incompleteness Theorems, we cannot disregard the consequent deeply inherent limitation of all axiomatic mathematical theories. \\
And that limitation cannot be kept away from nonrigid mathematical structures either, since such structures are also built as axiomatic theories. \\

And that G\"{o}delian limitation comes to further suggest the answer to the issue of what is Topology, is it the HKB one,  or is it one or another pseudo-topology ? \\
And the answer is simple indeed : the rigid HKB concept of topology may be just as little unique, as that of Geometry proved to be two centuries earlier ... \\

Regarding Mathematics in general, fortunately, two possible further developments, away from that G\"{o}delian limitation of all axiomatic theories, have recently appeared, even if they are not yet clearly enough in the general mathematical awareness. Namely, {\it self-referential} axiomatic mathematical theories, and perhaps even more surprisingly, {\it inconsistent} axiomatic mathematical theories, [14]. \\ \\

{\bf 1. Rings Are Non-Rigid Structures} \\

In Group Theory, given any group $( G, \diamond )$, be it commutative or not, all the concepts of the theory will in the last analysis be uniquely defined by the set $G$ and the binary operation $\diamond$. \\

On the other hand, in Ring Theory, given a ring $( R, +, .)$, the binary operation $''+''$ of addition, does {\it not} in general determine uniquely the binary operation $''.''$ of multiplication. Instead, they are only supposed to satisfy certain {\it compatibility} relations, namely, the distributivity of multiplication with respect to
addition. \\

And rather simple examples show that, given a commutative group, there can be more than one ring multiplication on
it. \\

Indeed, let $\mathbb{M}\,^n ( \mathbb{R} )$, with $n \geq 2$, be the set of $n \times n$ matrices of real numbers, and consider on it the commutative group structure given by the usual addition $''+''$ of matrices. \\

The following two {\it different} ring structures can be defined on $\mathbb{M}\,^n ( \mathbb{R} )$. \\

First, let $( \mathbb{M}\,^n ( \mathbb{R} ), +, . )$, where $''.''$ is the usual noncommutative multiplication of square matrices. Second, let $( \mathbb{M}\,^n ( \mathbb{R} ), +, \ast )$ where $''\ast''$ is the term by term multiplication of matrices, namely, given the matrices \\

$~~~~ A = ( a_{i, \, j} ~|~ 1 \leq i, j \leq n ),~~ B = ( b_{i, \, j} ~|~ 1 \leq i, j \leq n ) $ \\

we have $A . B = C$, where \\

$~~~~ C = ( c_{i, \, j} ~|~ 1 \leq i, j \leq n ) $ \\

with \\

$~~~~ c_{i, \, j} = a_{i, \, j} . b_{i, \, j} $ \\

And these two ring structures are indeed different, although their underlying commutative group structure is the same. For instance, the first one is noncommutative, while the second one is commutative. Furthermore, the first one has as unit element the square matrix with the diagonal 1, and with all the other elements 0, while the unit element in the second one is the matrix with all the elements 1. \\

Consequently, Ring Theory is indeed {\it non-rigid}. \\ \\

{\bf 2. The Ring $\mathbb{Z}$ Is Nearly Rigid} \\

Let $''+''$ denote the usual addition on $\mathbb{Z}$ while $''.''$ denotes the usual multiplication on it. Further, for a given integer \\

(2.1)~~~~ $ a \in \mathbb{Z} $ \\

let us consider on $\mathbb{Z}$ the binary operation $''\ast''$ defined by \\

(2.2)~~~~ $ n \ast m = a . n . m,~~~~ n, m \in \mathbb{Z} $ \\

{\bf Lemma 2.1.} \\

$( \mathbb{Z}, +, \ast )$ is a {\it commutative ring}. \\

{\bf Proof.} \\

We have for $n, m, k \in \mathbb{Z}$ the relations \\

$~~~~ n \ast ( m \ast k) = a . n . ( m \ast k ) = a . n . ( a . m . k ) = a . a . n . m . k $ \\

while \\

$~~~~ ( n \ast m ) \ast k = a . ( n \ast m ) . k = a . ( a . n . m ) . k = a . a . n . m . k $ \\

thus the associativity of $''\ast''$. Also we have the relations \\

$~~~~ n \ast ( m + k ) = a . n . ( m + k ) = ( a . n . m ) + ( a . n . k ) = ( n \ast m ) + ( n \ast k ) $ \\

hence the distributivity of $''\ast''$ with respect to $''+''$.

\hfill $\Box$ \\

Obviously, if $a = 1$, then $( \mathbb{Z}, +, \ast )$ is the usual ring $( \mathbb{Z}, +, . )$. On the other hand, if $a = -1$, then $( \mathbb{Z}, +, \ast )$ has the somewhat surprising multiplication rule \\

(2.3)~~~~ $ n \ast m = - n . m,~~~~ n, m \in \mathbb{Z} $ \\

and we shall call this the {\it alternate} ring of the usual ring $( \mathbb{Z}, +, . )$. \\

{\bf Lemma 2.2.} \\

$( \mathbb{Z}, +, \ast )$ is a {\it unital} ring, if and only if $a = \pm 1$, in which case it reduces to the usual ring  $( \mathbb{Z}, +, . )$, or to its alternate, see (2.3). \\

{\bf Proof.} \\

Let $u \in \mathbb{Z}$, such that \\

$~~~~ u \ast n = n,~ n \in \mathbb{Z} $ \\

then \\

$~~~~ a . u . n = n,~ n \in \mathbb{Z} $ \\

thus in particular \\

$~~~~ a . u = 1 $ \\

which, in view of (2.1), means that \\

$~~~~ a = u = 1 $, ~or~ $ a = u = -1 $

\hfill $\Box$ \\

Recalling now (2.1), (2.2), we obtain \\

{\bf Theorem 2.1.} \\

All the commutative ring structures on the commutative group $( \mathbb{Z}, + )$ are of the form $( \mathbb{Z}, +, \ast )$, for suitable $a \in \mathbb{Z}$. \\

{\bf Proof.} \\

Let be given any commutative ring  $( \mathbb{Z}, +, \circ )$, then we denote \\

(2.4)~~~~ $ a = 1 \circ 1 $ \\

Let now $n, m \in \mathbb{Z},~ n, m > 0$, then \\

(2.5)~~~~ $ n \circ m = a . n . m $ \\

Indeed \\

$~~~~ \begin{array}{l}
              1 \circ m = 1 \circ ( 1 + \ldots + 1 ) = a + \ldots + a = a . m \\ \\
              n \circ m = ( 1 + \ldots + 1 ) \circ m =
                           1 \circ m + \ldots + 1 \circ m = a . m + \ldots + a . m = a . n . m
      \end{array} $ \\ \\

{\bf Remark 2.1.} \\

1. As noted in section 1, Ring Theory is non-rigid, since the commutative group structure of a ring does not in general determine uniquely the ring multiplication. \\

However, as seen in Theorem 2.1. above, the usual commutative group structure on $\mathbb{Z}$ does determine the commutative multiplication on it, except for a constant factor in (2.1), (2.2). \\

As also seen in Lemma 2.1. above, the usual commutative group structure on $\mathbb{Z}$ does further determine the commutative multiplication on it in case this multiplication has a unit element, except for the possibility of the alternate ring structure. \\

2. The fact that on such a small set like $\mathbb{Z}$, small in the sense of having the smallest infinite cardinal, there are not many significantly different ring structures on its usual commutative group need not come as a surprise. Indeed, in [9,11] it was shown that on $\mathbb{N}$ there are few associative binary operations which satisfy some rather natural and mild conditions. \\

On the other hand, if we consider the commutative group $( \mathbb{M}\,^n ( \mathbb{Q} ), + )$ of $n \times n$ matrices of rational numbers with the usual addition of matrices, then as seen in section 1, two different ring structures can be associated with that group. Yet the set $\mathbb{M}\,^n ( \mathbb{Q} )$ has also the smallest infinite cardinal. It may therefore be the case that in the mentioned result in [9,11] the usual linear order on $\mathbb{Z}$, and thus induced on $\mathbb{N}$ as well, an order missing on $\mathbb{M}\,^n ( \mathbb{Q} )$, plays a role. \\

3. Obviously, the construction in (2.1), (2.2) can be applied to an arbitrary ring, and the result in Lemma 2.1. will still hold. \\

In particular, the {\it alternate} in (2.3) can be defined far arbitrary rings. \\

As for the result in Lemma 2.2., its proof uses the fact that in a unital ring $( R, +, . )$, the equation \\

(2.6)~~~~ $ a . u = 1,~~~~ a, u \in R $ \\

should only have the solutions \\

(2.7)~~~~ $ a = u = \pm 1 $ \\

Thus we have in general \\

{\bf Lemma 2.3.} \\

Given a unital ring $( R, +, . )$ with property (2.6), (2.7). If $a \in R$, then the ring $( R, +, \ast )$ obtained through (2.1), (2.2) is unital, if and only if $a = \pm 1$, thus it reduces to $( R, +, . )$, or to its alternate. \\

\end{document}